\documentclass[11pt,twoside]{amsart}

\author{Valery Alexeev} 
\address{Department of Mathematics\\
University of Georgia\\
Athens, GA 30602}
\email{valery@math.uga.edu}

\newcommand{\isoto}{{\overset{\sim}{\rightarrow}}}

\newcommand{\bP}{{\mathbb P}}
\newcommand{\bQ}{{\mathbb Q}}
\newcommand{\bZ}{{\mathbb Z}}

\newcommand{\bR}{{\mathbb R}}

\newcommand{\bG}{{\mathbb G}}

\newcommand{\cM}{{\mathcal M}}
\newcommand{\cP}{{\mathcal P}}

\newcommand{\cO}{{\mathcal O}}

\newcommand{\wC}{\widetilde C}

\newcommand{\Spec}{\operatorname{Spec}}

\newcommand{\Pic}{\operatorname{Pic}}

\newcommand{\Aut}{\operatorname{Aut}}

\theoremstyle{plain}
\newtheorem{thm}{Theorem}[section]

\newtheorem{lem}[thm]{Lemma}
\newtheorem{cor}[thm]{Corollary}

\theoremstyle{definition}
\newtheorem{defn}[thm]{Definition}
\newtheorem{defnprop}[thm]{Definition-Proposition}

\newtheorem{notation}[thm]{Notation}

\newtheorem{saynum}[thm]{}

\newtheorem{example}{Example}

\newtheorem{rem}[thm]{Remark}

\newtheorem{ack}{Acknowledgments}   

\theoremstyle{remark}

\newenvironment{say}{}{}

\newcommand{\Jac}{\operatorname{Jac}}
\newcommand{\gr}{\operatorname{gr}}

\newcommand{\ul}{\underline\lambda}
\newcommand{\uo}{\underline\omega}
\newcommand{\ud}{\underline d}
\newcommand{\ue}{\underline e}
\newcommand{\tg}{\tilde g}

\newcommand{\XR}{X\otimes\bR}

\newsymbol\rightsquigarrow 1320

\makeatletter
\input epic.sty
\input eepic.sty
\makeatother

\begin{document}
\bibliographystyle{amsalpha+}
\title[Compactified jacobians]{Compactified Jacobians}
\date{August 15, 1996}

\maketitle

{\bf
\begin{center}
  \begin{tabular}{@{}l@{}}
    Preliminary version
  \end{tabular}
\end{center}
}

\tableofcontents

\setcounter{section}{-1}

\section{Introduction}
\label{sec:Introduction}

\begin{saynum}
  Let $C$ be a reduced projective curve over a field $k$ such that
  $C_{\bar k}$ has only nodes as singularities. The jacobian
  $\Pic^0\,C$ is a semiabelian variety over $k$ which parameterizes
  invertible sheaves on $C$ of degree 0 on each irreducible component.
  It need not be proper.  The problem of finding a good
  compactification for it goes back at least to the work of Igusa
  \cite{Igusa_SystemsOfJacobians} and the notes of Mumford and Mayer
  \cite{Mumford_BoundaryPoints,Mayer69}. For an irreducible curve the
  answer was given already by D'Souza in \cite{DSouza79}.  Altman,
  Kleiman and others extended this work to the families of irreducible
  curves with more general (for example, nonplanar) singularities in a
  series of works \cite{AltmanKleiman_CompJac,AltmanKleiman_CompJac2,
    AltmanIarrobinoKleiman_IrreducibilityCompJac,
    KleimanKleppe81,Rego_CompactifiedJacobian}, and more recently
  \cite{Soucaris94,Esteves95}.
\end{saynum}
     
\begin{saynum}
  In the case when $C$ is reducible the situation is more complicated.
  In a classical paper \cite{OdaSeshadri79} Oda and Seshadri
  constructed a family of compactified jacobians $\Jac_{\phi}$
  parameterized by an element $\phi$ of a certain real vector space.
  The construction is very general and covers a lot of cases.  At the
  same time it poses a question of giving a more natural definition
  for $\Jac_{\phi}$ and explaining where exactly the multitude of
  answers comes from. A related paper is \cite{Seshadri82}.
\end{saynum}

\begin{rem}
  It is important to note that the term ``compactified jacobian'' is a
  misnomer. Most of the varieties $\Jac$ discussed here do not
  naturally contain $\Pic^0(C)$. This becomes especially clear when
  one works over a nonclosed field or with families.  Instead, there
  is always an action of $\Pic^0(C)$ and $\Jac$ is stratified into
  locally closed subschemes so that every stratum is a homogeneous
  space over $\Pic^0(C)$ and the maximal-dimensional strata are
  principal homogeneous spaces.
  
  However, we will use the term since it is widely accepted.
\end{rem}
  
\begin{saynum}
  The work of Simpson \cite{Simpson94a} on the moduli of coherent
  sheaves on projective schemes implies, as a very special case of a
  much more general situation, a natural definition of the
  compactified jacobian $\Jac_{d,L}(C)$ which depends on an integer
  $d$, the degree, and on an ample invertible sheaf $L$ on $C$, the
  polarization.  The definition is functorial and therefore also works
  for families. It turns out that $\Jac_{d,L}(C)$ and Oda-Seshadri's
  $\Jac_{\phi}(C)$ coincide and there is a simple formula for $\phi$
  as a function of $d$ and $L$. An immediate corollary of this is that
  they are all reduced and Cohen-Macaulay schemes.
\end{saynum}

\begin{saynum}
  In the case when the curve $C$ is stable, we can further narrow the
  choices by using for $L$ the dualizing sheaf $\omega_C$. Then for
  every $d\in\bZ$ the schemes $\Jac_{d,\omega_C}/\Aut(C)$ can be put
  in a family over the moduli space $\overline{M}_g$, where $g$ is the
  arithmetical genus of the curve $C$. This is the result of the work
  \cite{Pandharipande94} of Pandharipande (he also considers sheaves
  of rank $\ge2$).  A yet another family
  $\overline{P}_d\to\overline{M}_g$ for $d\ge10(2g-2)$ was earlier
  constructed by Caporaso in \cite{Caporaso94} as the compactification
  of the universal jacobian. The interpretation of the fiber
  $\overline{P}_d(C)$ over $[C]\in\overline{M}_g$ is in terms of
  invertible sheaves on certain semistable curves that have $C$ as a
  stable model.  Pandharipande shows that Caporaso's construction of
  $\overline{P}_d$ is equivalent to his.
\end{saynum}

\begin{saynum}
  Another approach is to look at a one-parameter family of smooth
  curves $C_t$ degenerating to $C$ and try to find a limit of the
  family of ``jacobians'' $\Jac_d(C_t)=\Pic^d(C_t)$, perhaps after a
  finite ramified base change. In the complex analytic situation
  Namikawa \cite{Namikawa_ToroidalDegsAVs} constructed infinitely many
  toroidal degenerations of principally polarized abelian varieties
  that depend on polyhedral decompositions. We note a related work
  \cite{Kajiwara93} where the compactified jacobians corresponding to
  polyhedral decompositions appear in the context of log geometry
  (under the restriction that the irreducible components of the curve
  $C$ are nonsingular).
  
  Among various polyhedral decompositions Namikawa explicitly
  distinguished one called the Voronoi decomposition (and the Delaunay
  decomposition dual to it).
\end{saynum}

\begin{saynum}
  The degeneration corresponding to the Delaunay and Voronoi
  decompositions also appears in \cite{AlexeevNakamura96} as a result
  of the ``simplified Mumford's construction''. There it is shown that
  a family of principally polarized abelian varieties with theta
  divisors over spectrum of a complete DVR has the canonical limit
  (perhaps after making a finite ramified base change first). This
  limit was called a stable quasiabelian variety (SQAV), and when
  considered as a pair $(P,B)$ with the theta divisor -- a stable
  quasiabelian pair (SQAP). This poses a question of whether the SQAP
  which appears as the limit of jacobians is one of the compactified
  jacobians $\Jac_{\phi}$, $\Jac_{d,L}(C)$ above, and if yes, then
  which one.
\end{saynum}

\begin{saynum}
  As explained in \cite{Alexeev_CMAV}, an SQAV corresponding to a
  smooth curve $C$ coincides with $\Pic^{g-1}(C)$ and not with
  $\Pic^0(C)$. This gives a motivation to look at the case $d=g-1$
  more closely.
\end{saynum}

\begin{saynum}
  We show that precisely for one degree, $d=g-1$, the scheme
  $\Jac_{d,L}(C)$ does not depend on the polarization $L$, so one can
  simply write $\Jac_{g-1}$. For this reason we call it the {\em
    canonical compactified jacobian.\/} We show that $\Jac_{g-1}$
  possesses a natural ample sheaf with a natural section which we call
  the theta divisor $\Theta_C$. We give a very simple explicit
  combinatorial description of the stratification of $\Jac_{g-1}$ and
  the restrictions of $\Theta_C$ on each stratum. The description goes
  in terms of the orientations on complete subgraphs of the dual graph
  $\Gamma(C)$ and invertible sheaves on the partial normalizations of
  $C$ of multidegrees that correspond to these orientations.
\end{saynum}

\begin{saynum}
  By considering a degenerating family $C_t\rightsquigarrow C$ of
  curves and the corresponding degenerating family $\Pic^{g-1}(C_t)
  \rightsquigarrow \Jac_{g-1}(C)$ one can show that the canonical
  compactified jacobian is an SQAV (we do not include this argument
  here). Therefore, the functor associating to each stable curve $C$
  its canonical compactified jacobian $(\Jac_{g-1},\Theta)$ should
  define a map from the Deligne-Mumford compactification
  $\overline{M}_g$ to the complete moduli space of SQAVs if the latter
  exists.
\end{saynum}

\begin{saynum}
  Finally, an SQAV was defined in
  \cite{AlexeevNakamura96,Alexeev_CMAV} in terms of some explicit
  combinatorial data. We give the corresponding data for a curve $C$
  (about a half of this description can already be found in
  \cite[\S18]{Namikawa_NewComp2} and
  \cite[9.D]{Namikawa_ToroidalCompSiegel} where it is attributed to
  Mumford). Further, we explain how this second description is related
  to the previous one.
\end{saynum}

\begin{saynum}
  Although the main definitions and results including
  \ref{defn:Jac_functor} and \ref{thm_Simpson} hold over arbitrary
  base, for most of the paper we will be working over an algebraically
  closed field for simplicity.
\end{saynum}

\begin{ack}
  I would like to thank Profs. Yu.I. Manin, T. Oda and R. Smith for
  very useful conversations. It was Prof. Oda who kindly explained me
  the relation between his and Seshadri's parameter $\phi$ and the
  polarizations (in the case of degree 0). Of course, it is completely
  my fault if I didn't get it right.
\end{ack}


\section{Definition of $\Jac_{d,L}$}
\label{sec:Definition of oJd}

\begin{say}
  In this section I give the definition for compactified jacobians
  which I feel is the easiest and the most natural, and formulate the
  main existence theorem, due to Simpson. To take the bull by the
  horns, here it is:
\end{say}

\begin{defn}\label{defn:cj}
  For every integer $d$ and a polarization $L$ on $C$, the
  ``compactified jacobian'' $\Jac_{d,L}$ is the coarse moduli space of
  semistable w.r.t. $L$ admissible sheaves on $C$ of degree $d$ up to
  the $\gr$-equivalence.
\end{defn}

\begin{say}
  Let us now explain the terms used in this definition. The {\em
    polarization\/} $L$ is an ample invertible sheaf. By {\em
    admissible\/} (for our purposes) sheaf we mean a coherent
  $\cO_C$-module $F$
  \begin{enumerate}
  \item which is of rank 1, i.e. it is invertible on a dense open
    subset of $C$, and
  \item such that for every $x\in C$ the stalk $F_{x,C}$
  is of depth 1. Equivalently, any nonzero subsheaf of $F$ has support
  of dimension 1.
  \end{enumerate}
  
  The latter condition is what Seshadri \cite{Seshadri82} calls a
  depth 1 sheaf and what Simpson \cite{Simpson94a} calls a purely
  dimensional sheaf.
\end{say}
\smallskip

\begin{saynum}\label{saynum:admissible_sheaves}
  As well known (see f.e. \cite{Seshadri82}) admissible sheaves have
  a very simple description:
  \begin{enumerate}
  \item if $x$ is nonsingular, $F_{x,C}\simeq \cO_{x,C}$
  \item if $x$ is a node, $F_{x,C}$ is either $\cO_{x,C}$ or the
    maximal ideal $m_{x,C}$. In the latter case $F$ is isomorphic to
    $\pi(x)_*F(x)$ where $\pi(x)$ is a partial normalization of $C$ at
    $x$ and $F(x)=\pi(x)^*F/\text{torsion}$.
  \end{enumerate}
  
  Moreover, if we are interested in the depth 1 sheaves that have rank
  0 or 1 at every generic point, we need to add the sheaves such that 
  \begin{enumerate}\setcounter{enumi}{2}
  \item $F_{x,C}=0$
  \item if $x$ is a node lying on two irreducible components $C_1$ and
    $C_2$ with the inclusions $i_k:C_k\to C$, then
    $F_{x,C}=i_{k*}\cO_{C_k}$, $k=1$ or $2$.
  \end{enumerate}
\end{saynum}

\begin{lem}\label{lem:admissible_sheaves}
  For each admissible sheaf $F$ on $C$ denote by $\pi': C'\to C$ the
  partial normalization of $C$ at the nodes where $F$ is not
  invertible and by $F'=\pi{'}^*F/tors$. Then $F=\pi'_*F'$. Therefore,
  every admissible sheaf $F$ on $C$ can be identified with a unique
  invertible sheaf $F'$ on a unique partial normalization $C'$ of $C$.
\end{lem}
\begin{proof}
  Well known, see f.e. \cite{Seshadri82}.
\end{proof}

\begin{defn}\label{defn:degree}
  For an admissible sheaf $F$ on $C$ the degree is defined as
  \begin{displaymath}
    \deg F= \chi(F) -\chi(\cO_C) =\chi(F) +g-1,
  \end{displaymath}
  where $g$ is the arithmetical genus of $C$.
\end{defn}

\begin{rem}
  Note that $\deg F'$ if defined on $C'$ itself is $\deg F$ minus the
  number of nodes where $F$ is not invertible.
\end{rem}

\begin{defn}
  Let $\ul=(\lambda_1\dots\lambda_s)$ be the multidegree of $L$, and
  $\underline{r}=(r_1\dots r_s)$ be the multirank of a depth 1 sheaf
  $F$. The Seshadri slope is
  \begin{displaymath}
    \mu_L(F)= \frac{\chi(F)}{\sum \lambda_i r_i}
  \end{displaymath}
\end{defn}

\begin{defn}
  A depth 1 sheaf $F$ on $C$ is said to be stable (resp. semistable)
  w.r.t. the polarization $L$ if for any nonzero subsheaf $E\subset F$
  one has
  \begin{displaymath}
    \mu_L(E) < \mu_L(F)
  \end{displaymath}
  (resp $\le$).
\end{defn}

\begin{rem}
  This definition, due to Seshadri \cite{Seshadri82}, is a particular
  case of a much more general one given by Simpson in
  \cite{Simpson94a} that applies to a pure-dimensional sheaf on any
  projective scheme whatsoever.
\end{rem}

\begin{lem}\label{lem:criterion_stability_for_admissible_sheaves}
  A depth 1 sheaf is (semi)stable iff the inequality $\mu_L(E)<$(resp.
  $\le$) $\mu_L(F)$ is satisfied for finitely many subsheaves of the
  form $F_D=i_*(i^*F/tors)$ for every subcurve $D\subset C$, where
  $i:D\to C$ is the inclusion morphism.
\end{lem}
\begin{proof}
  Indeed, for any depth 1 subsheaf $E$ with support $D$ one has
  $E\subset F_D$ and $\mu_L(E)\le\mu_L(F_D)$.
\end{proof}

\begin{rem}
  This leads to a series of simple inequalities some of which will be
  considered in the next sections. Therefore, knowing the multidegree
  of $F$ and the set of nodes where $F$ is not locally free, it is
  easy to say whether $F$ is (semi)stable or not.
\end{rem}

\begin{say}
  According to the general theory, for every depth 1 sheaf $F$ there
  exists a Harder-Narasimhan filtration
  \begin{displaymath}
    0=F_0\subset F_1 \subset\dots\subset F_k=F
  \end{displaymath}
  with strictly decreasing slopes and semistable quotients
  $F_i/F_{i+1}$.
  
  If $F$ is semistable, then there is a similar (Jordan-Holder)
  filtration with stable $F_i/F_{i+1}$ which is not unique. However,
  the graded object $\gr(F)=\oplus F_i/F_{i+1}$ is uniquely defined.
\end{say}

\begin{defn}
  Two semistable depth 1 sheaves $F$ and $F'$ are said to be
  $\gr${\em-equivalent\/} if $\gr(F)\simeq\gr(F')$.
\end{defn}

\begin{say}
  Now all the ingredients of the definition \ref{defn:cj} have been
  introduced.
\end{say}
\smallskip

\begin{say}
  To put this into a functorial perspective, consider a projective
  morphism of schemes $\pi:C\times S\to S$ whose every geometric fiber
  is a reduced curve with nodes only as singularities and a relatively
  ample sheaf $L$ on $C$. We will say that a coherent sheaf $F$ on $C$
  is admissible (resp. stable, resp.  semistable) if so is its
  restriction to every geometric fiber of $\pi$. We say that two
  sheaves are equivalent (resp.  $\gr$-equivalent) if the restrictions
  on the geometric fibers are isomorphic (resp.  $\gr$-equivalent).
  Now define the moduli functor $Jac_{d,L}(C/S):Schemes\to Sets$ in
  the following way:
\end{say}

\begin{defn}\label{defn:Jac_functor}
  For any scheme $S'$, $Jac^-_{d,L}(C/S)(S')$ is the set of semistable
  admissible sheaves on $C'=C\underset{S}{\times}S'/S'$ up to the
  $\gr$-equivalence. The functor $Jac^-_{d,L}(C/S)$ itself is not
  necessarily a sheaf for the fppf (faithfully flat of finite
  presentation) topology and therefore cannot be representable even if
  all the sheaves are stable. This happens basically for the same
  reason why the functor $\Pic\,S'$ is not a sheaf: one needs a
  rigidification to kill the infinite ($=\bG_m$) group of
  automorphisms.
  
  When the smooth locus of $C/S$ has a section, one can use the
  rigidified version. Or we can follow the same path as for the
  relative Picard functor, i.e. define $Jac_{d,L}(C/S)$ to be the
  fppf-sheafification of $Jac^-_{d,L}(C/S)$.
\end{defn}

\begin{thm}[Simpson]\label{thm_Simpson}
  The functor $Jac_{d,L}(C/S)$ is coarsely represented by a projective
  scheme $\Jac_{d,L}(C/S)$.
\end{thm}
\begin{proof}
  This may be proved by the same methods as in \cite{Simpson94a} and
  basically is a very special case of \cite[1.21]{Simpson94a}. Simpson
  works over $\mathbb C$ but in our particular situation there is no
  need for this. The hardest question involved, the boundedness, is
  basically obvious.
  
  Alternatively, in \ref{thm:Jacs_are_the_same} we will show that our
  (semi)stability condition is equivalent to the one that was used by
  Oda and Seshadri, so in the case $S=\Spec\bar k$ the theorem follows
  by \cite[12.14]{OdaSeshadri79}.
\end{proof}


\section{Basic definitions and notations}
\label{sec:Basic definitions and notations}

\begin{say}
  The purpose of this section is to fix some notations common to the
  following sections and to introduce the basic examples of curves on
  which the later descriptions will be illustrated.
\end{say}

\begin{notation}
  \begin{enumerate}
  \item To any curve $C$ we can associate the unoriented graph
    $\Gamma(C)$ by assigning a vertex to each irreducible component
    $C_i$ and an edge to every node. We do not assume $C$ to be
    connected, so the graph need not be connected either.
  \item $\pi:\widetilde C\to C$ denotes the normalization of $C$.
  \item $C_i$ are the irreducible components of $C$, $\widetilde C_i$
    are their normalizations.
  \item $g_i=p_a(C_i)$, $\tilde g_i=p_a(\widetilde C_i)$.
  \item $h(C)=h(\Gamma(C))$ is the cyclotomic number -- the number of
    independent loops in $\Gamma$, i.e. the rank of $H_1(\Gamma(C))$
    when $\Gamma$ is considered as a cell complex.
  \end{enumerate}
\end{notation}

{\bf Six simple examples.}
  
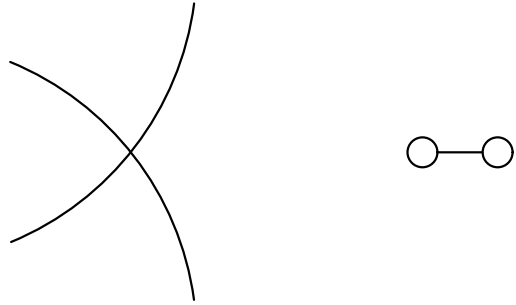
\setlength{\unitlength}{0.1cm}
\begin{picture}(110,55)(0,0)
\thicklines
\put(0,0){\arc{80}{5.1025}{6.1492}}
\put(0,50){\arc{80}{0.1308}{1.1775}}
\put(70,25){\circle{4}}
\put(80,25){\circle{4}}
\drawline(72,25)(78,25)
\end{picture}

\begin{example}\label{example1}
  $C=C_1\cup C_2$ intersecting at one point with both $C_i$
  smooth. $\tg_i=g_i$, $h(C)=0$. 
\end{example}

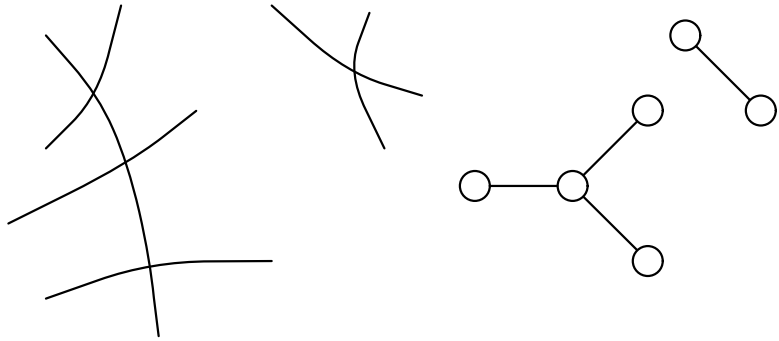
\setlength{\unitlength}{0.1cm}
\begin{picture}(110,55)(0,0)
\thicklines
\spline(10,45)(17,37)(21,28)(24,15)(25,5)
\spline(10,30)(17,37)(20,49)
\spline(5,20)(21,28)(30,35)
\spline(10,10)(24,15)(40,15)
\spline(40,49)(50,40)(60,37)
\spline(53,48)(50,40)(55,30)
\put(67,25){\circle{4}}
\put(80,25){\circle{4}}
\put(90,35){\circle{4}}
\put(90,15){\circle{4}}
\drawline(69,25)(78,25)
\drawline(81.5,26.5)(88.5,33.5)
\drawline(81.5,23.5)(88.5,16.5)
\put(95,45){\circle{4}}
\put(105,35){\circle{4}}
\drawline(96.5,43.5)(103.5,36.5)
\end{picture}
\begin{example}\label{example2}
  The generalization of the previous example is a curve whose dual
  graph is a forest.  Still $\tg_i=g_i$ and $h(C)=0$.
\end{example}

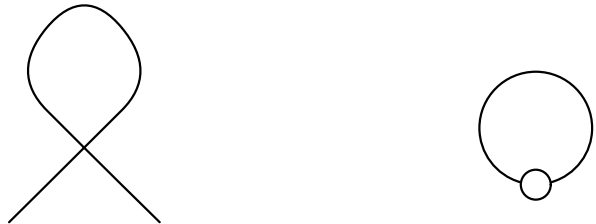
\setlength{\unitlength}{0.1cm}
\begin{picture}(110,55)(0,0)
\thicklines
\spline(20,10)(30,20)(40,30)(30,42)(20,30)(30,20)(40,10)
\put(90,15){\circle{4}}
\put(90,22.5){\arc{15}{1.8555}{7.5647}}
\end{picture}
\begin{example}\label{example3}
  An irreducible curve with one node. $\tg=g-1$, $h(C)=1$. 
\end{example}

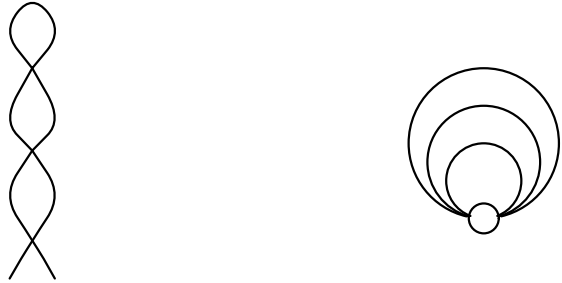
\setlength{\unitlength}{0.1cm}
\begin{picture}(110,55)(0,0)
\thicklines
\spline(30,35)(26,40)(30,45)(34,40)(30,35)
\spline(30,35)(34,28)(30,24)
\spline(30,35)(26,28)(30,24)
\spline(30,24)(34,18)(30,12)(27,7)
\spline(30,24)(26,18)(30,12)(33,7)
\put(90,15){\circle{4}}
\put(90,20){\arc{10}{1.9327}{7.4871}}
\put(90,22.5){\arc{15}{1.8454}{7.5744}}
\put(90,25){\arc{20}{1.8019}{7.6179}}
\end{picture}
\begin{example}\label{example4}
  The generalization of that is an irreducible curve with $n$
  nodes. $\tg=g-n$, $h(C)=n$. 
\end{example}

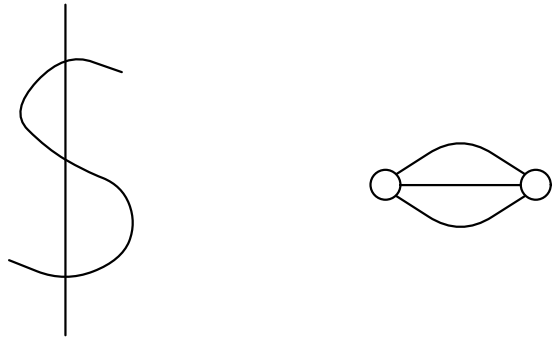
\setlength{\unitlength}{0.1cm}
\begin{picture}(110,55)(0,0)
\thicklines
\spline(20,15)(27.5,12)(35,15)(37,20)(35,25)(30,27)(25,30)(20,35)
(27.5,43)(35,40)
\drawline(27.5,49)(27.5,5)
\put(70,25){\circle{4}}
\put(90,25){\circle{4}}
\drawline(72,25)(88,25)
\spline(71.5,26.5)(80,32)(88.5,26.5)
\spline(71.5,23.5)(80,18)(88.5,23.5)
\end{picture}
\begin{example}\label{example5}
  The dollar sign curve. $\tg_i=g_i$, $h(C)=2$. 
\end{example}

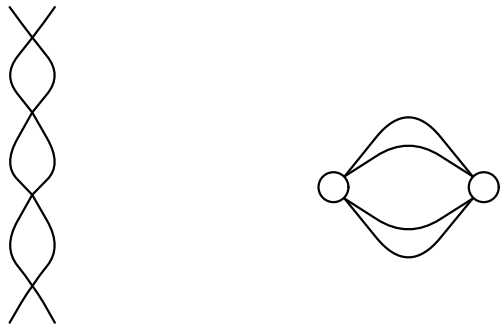
\setlength{\unitlength}{0.1cm}
\begin{picture}(110,55)(0,0)
\thicklines
\spline(30,35)(26,40)(30,45)(33,49)
\spline(30,35)(34,40)(30,45)(27,49)
\spline(30,35)(34,28)(30,24)
\spline(30,35)(26,28)(30,24)
\spline(30,24)(34,17)(30,12)(27,7)
\spline(30,24)(26,17)(30,12)(33,7)
\put(70,25){\circle{4}}
\put(90,25){\circle{4}}
\spline(71.5,26.5)(80,32)(88.5,26.5)
\spline(71.5,23.5)(80,18)(88.5,23.5)
\spline(71.5,26.5)(80,37)(88.5,26.5)
\spline(71.5,23.5)(80,13)(88.5,23.5)
\end{picture}
\begin{example}\label{example6}
  The generalization of the dollar sign curve is a curve $C=C_1\cup
  C_2$ with both $C_i$ smooth and intersecting at $n$ points.
  $\tg_i=g_i$, $h(C)=n-1$.
\end{example}


\section{Comparison with Oda-Seshadri's compactified jacobians}
\label{sec:Comparison with Oda-Seshadri's compactified jacobians}

\begin{say}
  We would like to compare the compactified jacobians introduced in
  section \ref{sec:Definition of oJd} with those appearing in the
  classical paper \cite{OdaSeshadri79} of Oda and Seshadri.
\end{say}

\begin{saynum}
  The $\Jac_{\phi}$ in \cite{OdaSeshadri79} are constructed using GIT
  as the moduli spaces of $\phi$-semistable admissible sheaves.  Here
  $\phi$ is an element of a certain real vector space $\partial
  C_1(\Gamma,\bR)$ (without loss of generality one can assume that
  $\phi\in\partial C_1(\Gamma,\bQ)$). $\Gamma$ is, as in the previous
  section, the dual graph of $C$ and $C_0,C_1,H_0,H_1,C^0,C^1,H^0$ and
  $H^1$ are the associated to it chain and (co)homology groups.
\end{saynum}

\begin{saynum}
  In \cite{OdaSeshadri79} the main object of interest is the depth 1
  sheaves of degree 0. Oda and Seshadri give the combinatorial
  definition of a $\phi$-stable (resp. semistable) sheaf and introduce
  the $\phi$-equivalence relation. The main result then is that there
  exists a reduced scheme $\Jac_{\phi}$ which coarsely represents the
  functor of $\phi$-semistable sheaves up to $\phi$-equivalence.  This
  is then applied to compactify $\Pic^0\,C$.
  
  However, for any depth 1 sheaf of arbitrary degree $d$ one can
  relate the $\phi$-(semi)stability and equivalence with
  $(d,L)$-(semi)stability and equivalence. Then whatever is proved for
  $\Jac_{\phi}$ immediately applies to $\Jac_{d,L}$. Here is the
  precise connection.
\end{saynum}
 
\begin{thm}\label{thm:Jacs_are_the_same}
  Let $\ul=(\lambda_i)$ and $\uo=(\lambda_i)$ be the multidegrees of
  the polarization $L$ and of the dualizing sheaf $\omega_C$
  respectively, and $\lambda=\sum\lambda_i$ and
  $\omega=\sum\omega_i=2g-2$ be the total degrees.  Pick arbitrary
  integers $d_i$ with $\sum d_i=d$ and sufficiently large integers
  $\tilde{n}_i$.  Define $\phi=(\phi_i)\in\partial C_1(\bQ)$ to be a
  solution of the following system of linear equations
  \begin{displaymath}
    (\lambda_i/\lambda)(d-\omega/2)=
    d_i-\omega_i/2+\tilde{n}_i+\phi_i
  \end{displaymath}
  ($\phi$ is only defined up to a shift by a lattice).  Then an
  admissible sheaf of degree $d$ is (semi)stable w.r.t.  $L$ iff it is
  $\phi$-(semi)stable. Two semistable w.r.t. $L$ sheaves are
  $\gr$-equivalent iff they are $\phi$-equivalent.
\end{thm}
\begin{proof}
  This can be extracted from \cite[\S11]{OdaSeshadri79} directly,
  particularly from the account on pp.52-53.
\end{proof}

\begin{cor}
  Every $\Jac_{d,L}$ is isomorphic to one of $\Jac_{\phi}$ and vice
  versa. 
\end{cor}

\begin{cor}
  Every $\Jac_{d,L}$ is reduced and Cohen-Macaulay.
\end{cor}
\begin{proof}
  Indeed, $\Jac_{\phi}$ is reduced by \cite[11.4]{OdaSeshadri79}.
  Moreover, the proof shows (pp.60-62) that $\Jac_{\phi}$ is a good
  GIT quotient of a certain scheme $R$ and there exists an open
  subscheme $Y\subset R\times \bP(E^*)$ such that the projection $R\to
  Y$ is surjective, and $Y$ is formally smooth over a Hilbert scheme
  $H$ which is open in a quotient by the symmetric group of
  $C\times\dots\times C$.

  Therefore, $H$ is CM, an so is $Y$, and so is $R$, and so is
  $\Jac_{\phi}$. 
\end{proof}


\section{Description of $\Jac_{g-1}$}
\label{sec:first description}

\begin{lem}
  $\Jac_{g-1,L}(C)$ does not depend on the polarization $L$.
\end{lem}
\begin{proof}
  Indeed, by definition \ref{defn:degree} the degree $d=g-1$ iff
  $\chi(F)=0$. Then for any $E\subset F$ the inequality
  \begin{displaymath}
    \mu_L(E)\le \,(\text{resp. }<) \,\mu_L(F)
  \end{displaymath}
  is equivalent to
  \begin{displaymath}
    \chi(E)\le \,(\text{resp. }<) \,0
  \end{displaymath}
\end{proof}

\begin{say}
  For this reason we will call $\Jac_{g-1}(C)$ the {\em canonical
    compactified jacobian.}
\end{say}

\begin{defn}
  A subgraph $\Gamma'\subset\Gamma$ is said to be {\em generating\/}
  if
  $\operatorname{vertices}(\Gamma)=\operatorname{vertices}(\Gamma')$.
  Every such subgraph corresponds to a partial normalization of $C$ at
  the nodes $\Gamma-\Gamma'$. We denote this partial normalization by
  $\pi(\Gamma'):C(\Gamma')\to C$. Note in particular that
  $C(\Gamma)=C$ and that $\tilde C$ is $C(\Gamma')$ where $\Gamma'$
  has all the vertices of $\Gamma$ but no edges at all.
\end{defn}

\begin{defn}
  A subgraph $\Gamma'\subset\Gamma$ is said to be {\em complete\/} if
  $\operatorname{vertices}(\Gamma')\subset
  \operatorname{vertices}(\Gamma)$ and $\operatorname{edges}(\Gamma')$
  are precisely the edges of $\Gamma$ lying inside $\Gamma'$. These
  graphs correspond to subcurves $D\subset C$. Often we identify such
  subcurves $D$ with the corresponding subgraphs.
\end{defn}

\begin{defn}
  A {\em multidegree\/} of a graph $\Gamma$ is a set
  $\underline{d}=(d_i)$ of integers for every vertex $C_i$ of
  $\Gamma$. {\em We will always assume that\/}
  \begin{displaymath}
    \sum d_i = g-1
  \end{displaymath}
  
  A {\em normalization of multidegree\/} $\underline{d}$ is a set of
  integers $\underline{e}=(e_i)$ defined by
  \begin{displaymath}
    e_i=d_i-(\tilde g_i-1).
  \end{displaymath}
  It will be called {\em the normalized multidegree.\/} Note that we can
  use multidegrees $\underline{d}$ and normalized multidegrees
  $\underline{e}$ interchangeably. Note that $\sum e_i$ equals the
  number of edges of $\Gamma$.

  For a subcurve $D\subset C$, i.e. a complete subgraph
  $\Gamma'\subset\Gamma$, we set 
  \begin{displaymath}
    d_D=\sum_{C_i\subset D}  d_i, \quad
    e_D=\sum_{C_i\subset D}  e_i
  \end{displaymath}
\end{defn}

\begin{defnprop}\label{defnprop}
  A normalized  multidegree $\underline{e}$ is called {\em
    semistable\/} (resp. {\em stable\/}) if any of the following
  equivalent conditions hold:
  \begin{enumerate}
  \item 
    \begin{displaymath}
      |e_D-\#edges(D)-\frac{1}{2}D(C-D)|\le \frac{1}{2}D(C-D)
    \end{displaymath}
    for every subcurve $D\subset C$. Here $D(C-D)$ is the number of
    points in $D\cap\overline{(C\setminus D)}$ (resp. $<$).
  \item 
    \begin{displaymath}
      e_D\le \#edges(D)+D(C-D)
    \end{displaymath}
    (resp.$<$).
  \item there exists an orientation of the graph $\Gamma$ such that
    $e_i$ equals the number of arrows pointing at $C_i$ (resp. in
    addition there is no proper subcurve $D\subset C$ such that all
    arrows between $D$ and $C-D$ go in one direction).
  \end{enumerate}
  In this case the  multidegree $(\ud)$ is also called
  (semi)stable. 
\end{defnprop}
\begin{proof}
  The implication (i)$\Rightarrow$(ii) is clear and
  the inverse is obtained by looking at $D'=C-D$. (iii)
  obviously implies (ii).
  
  To prove the implication (ii)$\Rightarrow$(iii) first assume that
  the normalized multidegree $\underline{e}$ of the graph of $C$ is
  strictly semistable, i.e. there exists a subcurve $D\subset C$ for
  which the equality holds. Then consider separately the following
  multidegrees on $D$ and $C-D$. On $C-D$ simply take the restriction of
  $\underline{e}$. On $D$, however, for every vertex $C_i$ take
  $e_i'=e_i$ minus the number of edges between $C_i$ and $C-D$. Then
  it is easy to show that the two multidegrees thus obtained are
   semistable. Therefore, the orientations on $D$ and
  $C-D$ exist by the induction on the number of vertices. To complete
  the orientation of $C$, orient all the edges between $D$ and $C-D$
  to point at $D$.
  
  In general, starting with a semistable multidegree as in (ii) we can
  fix an arbitrary vertex $C_{i_0}$ and change the degrees of
  $C_{i_0}$ and the neighboring vertices = curves $C_j$ by 1 to make
  the multidegree strictly semistable, thus reducing to the previous
  case. Hence, we get an orientation for the modified multidegree.
  The orientation for the original multidegree is then obtained by
  reversing the orientations of edges $(i_0,j)$.
\end{proof}

\begin{say}
  The third condition of the above definition is the easiest to
  check. We will call an orientation satisfying (iii)
  semistable (resp. stable). Note that different orientations may well
  produce the same multidegree.
\end{say}
\smallskip

\begin{saynum}
  This is how the above combinatorial definitions relate to the
  (semi) stability of admissible sheaves on $C$. By lemma
  \ref{lem:admissible_sheaves} every admissible sheaf $F$ on $C$ can
  be identified with a unique invertible sheaf $F'$ on a unique
  partial normalization $C'=C(\Gamma')$ of $C$. Denote by $(\ud')$
  (resp.  $(\ue')$) the corresponding (resp. normalized) multidegrees on
  $C'$.  Then
\end{saynum}

\begin{lem}
  If $\deg F=g-1$, then for $(\ud')$ one has $\sum d'_i=g'-1$.
\end{lem}
\begin{proof}
  Obvious.
\end{proof}

\begin{lem}
  \begin{enumerate}
  \item $F$ is semistable iff $(\ue')$ is semistable.
  \item $F$ is stable iff $(\ue')$ is stable and the graphs $\Gamma$
    and $\Gamma'$ have the same number of connected components.
  \end{enumerate}
\end{lem}
\begin{proof}
  Follows easily from \ref{defnprop} and
  \ref{lem:criterion_stability_for_admissible_sheaves}. 
\end{proof}

\begin{say}
  We can now describe the points of $\Jac_{g-1}(C)$ as follows
\end{say}

\begin{thm}\label{thm:1description}
  \begin{enumerate}
  \item $\Jac_{g-1}(C)$ has a natural stratification into homogeneous
    spaces over $\Pic^0(C)$.  Each stratum corresponds in a 1-to-1 way
    to a stable multidegree $\underline{d'}$ (resp. stable normalized
    multidegree $\underline{e'}$) on a generating subgraph
    $\Gamma'\subset\Gamma$.  The $k$-points of this stratum can be
    identified with $k$-points of $\Pic_{\underline{d'}}(C(\Gamma'))$,
    i.e. with invertible sheaves on $C(\Gamma')$ of multidegree
    $\underline{d'}$. The codimension of this stratum equals
    $h(\Gamma)-h(\Gamma')$.
  \item There is a natural Cartier divisor $\Theta$ on
    $\Jac_{g-1}(C)$.  Under the above identification, the restriction
    of $\Theta$ on each stratum corresponds to the sheaves $L$ with
    $h^0(L)>0$.
  \end{enumerate}
\end{thm}

\begin{say}
  To illustrate this theorem, let us see what happens in our basic
  examples. 
\end{say}

\setcounter{example}{0}

\begin{example}
  Graph $\Gamma$ doesn't have any stable multidegrees: take $D$ to be
  one of the vertices. The only possibility then is $\Gamma'$ which is
  a disjoint union of two vertices and the multidegree
  $\underline{e}=(0,0)$, i.e.  $\underline{d}=(\tilde g_1-1,\tilde
  g_2-1)=(g_1-1,g_2-1)$. The graph $\Gamma'$ corresponds to the
  normalization $\widetilde X=X_1\bigsqcup X_2$ and
    \begin{displaymath}
      \Jac_{g-1}(C)=\Pic^{g_1-1}(C_1)\oplus\Pic^{g_2-1}(C_2)
    \end{displaymath}
\end{example}

\begin{example}
  Once again, a forest doesn't have any stable orientations unless all
  the brunches, i.e. edges, are cut. So, there is only one normalized
  multidegree $\ue'=(0,\dots,0)$ for a generating subgraph
  corresponding to the normalization $\widetilde C$ and
  \begin{displaymath}
    \Jac_{g-1}(C)=\oplus_i\Pic^{g_i-1}(C_i)
  \end{displaymath}
\end{example}

\begin{example}
  The stable orientations are 

  \setlength{\unitlength}{0.1cm}
  \begin{picture}(110,55)(0,0)
    \thicklines
    \put(15,15){\circle{4}}
    \put(15,22.5){\arc{15}{1.8555}{7.5647}}
    \drawline(16.75,15.25)(19,17.5)
    \drawline(16.75,15.5)(19.5,15.5)
    \put(35,15){\circle{4}}
  \end{picture}
  
  The first corresponds to $\Pic^{\tg-1+1}(C)=\Pic^{\tg}(C)$ and the
  second -- to $\Pic^{\tg-1}(\widetilde X)$.
\end{example}

\begin{example}
  There are $2^n$ subgraphs $\Gamma'$: each edge is either included in
  $\Gamma'$ or it's not. Each graph with $k$ edges obviously defines
  the multidegree $(\ud')=(k)$. Therefore there are $\binom n{n-k}=\binom
  nk$ strata of codimension $n-k$ each corresponding to
  $\Pic^{g'-1}(C')$.
\end{example}

\begin{example}[Dollar sign curve]
  The possible generating subgraphs are:

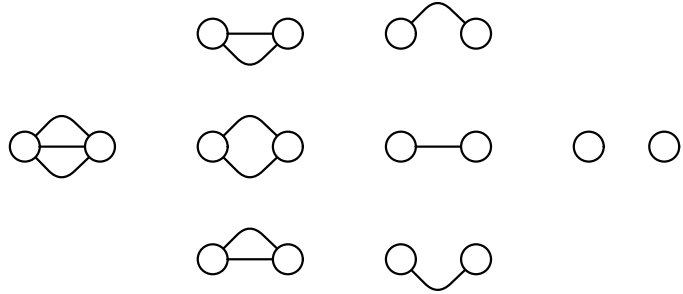
\setlength{\unitlength}{0.1cm}
\begin{picture}(110,55)(0,0)
\thicklines
\put(15,25){\circle{4}}
\put(25,25){\circle{4}}
\drawline(17,25)(23,25)
\spline(16.5,26.5)(20,30)(23.5,26.5)
\spline(16.5,23.5)(20,20)(23.5,23.5)
\put(40,25){\circle{4}}
\put(50,25){\circle{4}}
\spline(41.5,26.5)(45,30)(48.5,26.5)
\spline(41.5,23.5)(45,20)(48.5,23.5)
\put(40,40){\circle{4}}
\put(50,40){\circle{4}}
\drawline(42,40)(48,40)
\spline(41.5,38.5)(45,35)(48.5,38.5)
\put(40,10){\circle{4}}
\put(50,10){\circle{4}}
\drawline(42,10)(48,10)
\spline(41.5,11.5)(45,15)(48.5,11.5)
\put(65,40){\circle{4}}
\put(75,40){\circle{4}}
\spline(66.5,41.5)(70,45)(73.5,41.5)
\put(65,25){\circle{4}}
\put(75,25){\circle{4}}
\drawline(67,25)(73,25)
\put(65,10){\circle{4}}
\put(75,10){\circle{4}}
\spline(66.5,8.5)(70,5)(73.5,8.5)
\put(90,25){\circle{4}}
\put(100,25){\circle{4}}
\end{picture}

It is very easy to list all stable orientations and the corresponding
stable multidegrees. Here are some of them:

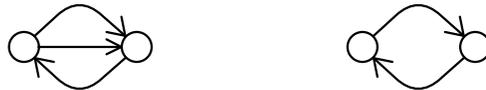
\setlength{\unitlength}{0.1cm}
\begin{picture}(110,55)(0,0)
\thicklines
\put(25,25){\circle{4}}
\put(40,25){\circle{4}}
\drawline(27,25)(38,25)
\drawline(38,25)(36,26)
\drawline(38,25)(36,24)
\spline(26.5,26.5)(32.5,32)(38.5,26.5)
\drawline(38.5,26.5)(38,29)
\drawline(38.5,26.5)(36,27)
\spline(26.5,23.5)(32.5,18)(38.5,23.5)
\drawline(26.5,23.5)(27,21)
\drawline(26.5,23.5)(29,23)
\put(70,25){\circle{4}}
\put(85,25){\circle{4}}
\spline(71.5,26.5)(77.5,32)(83.5,26.5)
\drawline(83.5,26.5)(83,29)
\drawline(83.5,26.5)(81,27)
\spline(71.5,23.5)(77.5,18)(83.5,23.5)
\drawline(71.5,23.5)(72,21)
\drawline(71.5,23.5)(74,23)
\end{picture}

Here is the complete list:
\begin{enumerate}
\item For the graph $\Gamma$ itself there are two multidegrees $(2,1)$
  and $(1,2)$ corresponding to invertible sheaves of multidegree
  $(\tg_1+1,\tg_2)=(\tg-1,\tg_2-1)+(2,1)$ and
  $(\tg_1,\tg_2+1)=(\tg-1,\tg_2-1)+(1,2)$ on $C$.
\item In the second column, for each graph we have a unique stable
  multidegree $(1,1)$. Hence, there are 3 strata of codimension 1
  corresponding to invertible sheaves of multidegree
  $(\tg_1,\tg_2)=(\tg-1,\tg_2-1)+(1,1)$.
\item In the third column there are no stable orientations -- the
  graphs are trees. 
\item From the last column we get the normalized multidegree $(0,0)$
  which corresponds to the invertible sheaves on the normalization
  $\widetilde C$ of $C$ of multidegree
  $(\tg_1-1,\tg_2-1)=(\tg-1,\tg_2-1)+(0,0)$.
\end{enumerate}

\end{example}

\begin{example}
  This is an exercise no harder then the previous five.  Any subgraph
  $\Gamma'$ has at least one stable orientation with one exception:
  when $\Gamma'$ contains only one edge. For each such subgraph with
  $k$ edges the number of possible stable multidegrees is $k-1$.
  Therefore, there are $\binom nk(k-1)$ strata of codimension $n-k$
  for $k>0$ and one stratum for $k=0$.
\end{example}

\begin{proof}[Proof of \ref{thm:1description}]
  The proof of (i) follows immediately from
  \ref{lem:admissible_sheaves} and
  \ref{lem:criterion_stability_for_admissible_sheaves}. In each
  $\gr$-equivalence class of strictly semistable sheaves we can choose
  the one with the minimal graph $\Gamma'$ and it will be stable in
  our definition.
  
  Next, we have to show the existence of a natural line bundle with a
  natural section on $\Jac_{g-1}(C)$.
  
  To define the divisor $\Theta$ in a way similar to how it was done
  in \cite{Soucaris94,Esteves95} for irreducible $C$.  Consider any
  family $\pi:C\times S\to S$ and an admissible sheaf $F$ of degree
  $g-1$ on $C\times S$. Then there is a natural line bundle $L(F)$ on
  $S$ defined as
  \begin{displaymath}
    L(F)= (\det R\pi_*F)^{-1}
  \end{displaymath}
  see \cite{KnudsenMumford76} for its definition.
  
  If we replace $F$ by $F\otimes\pi^*E$, $L(F)$ will be replaced by
  $L(F)\otimes E^{-\chi(F_t)}$. When $d=g-1$, $\chi=0$ which means
  that $L(F)$ will not change, so it is universally defined. Moreover,
  two $\gr$-equivalent families of semistable sheaves produce the same
  $L(F)$. The latter follows from the fact that if
  \begin{displaymath}
    0\to F'\to F\to F''\to 0
  \end{displaymath}
  is an exact sequence, then $\det R\pi_*F=(\det R\pi_*F')\otimes
  (\det R\pi_*F'')$, so only the stable factors are important.
  
  $\Jac_{d,L}$ and $\Jac_{\phi}$ are constructed using GIT as a
  quotient of the Grothendieck's $Quot$-schemes. By the universality
  $L(F)$ descends to $\Jac_{g-1}$.
  
  Now fix a point $c\in C$ and consider a universal family $F$ of
  invertible sheaves of degree $g'-1$ and multidegree $\ud'$ over
  $\Pic^{g'-1}(C')$, where $C'=C(\Gamma')$ is any of the partial
  normalizations of $C$. When does the formula
  \begin{displaymath}
    \Theta=\{s\in\Pic^{g'-1}(C') \,|\, h^0(F_s)>0\}
  \end{displaymath}
  define a divisor? The answer is given by a theorem of Beauville
  \cite[2.1]{Beauville_PrymSchottky}: it is exactly when the
  multidegree $\ud'$ is semistable using the part (iii) of
  Definition-Proposition \ref{defnprop}. It is also easy to show
  directly that if two sheaves of degree $g'-1$ are semistable and
  $\gr$-equivalent, then $h^0(F_1)\ne0$ iff $h^0(F_2)\ne0$.
  
  $\Theta$ provides a section of $(\det R\pi_*F)^{-1}$.
\end{proof}

\begin{rem}
  From the above proof we have a yet another characterization of
  semistable admissible sheaves in degree $g-1$: they have the
  multidegrees for which the usual definition of the theta-divisor
  actually gives a divisor.
\end{rem}


\section{An SQAV corresponding to a curve}
\label{sec:Jac second description}

\begin{saynum}
  An SQAV was defined in \cite{AlexeevNakamura96} explicitly starting
  from the following combinatorial data:
  \begin{enumerate}
  \item a lattice $X\simeq\bZ^{g'}$ (and a lattice $Y$ isomorphic to
    it via $\phi:Y\isoto X$).
  \item a symmetric positive definite bilinear form $B:X\times
    X\to\bZ$. 
  \item an abelian variety $A_0$ of dimension $g''$, $g'+g''=g$, with
    a principal polarization given by an ample sheaf $\cM_0$.
  \item a homomorphism $c_0:X\to A_0^t(k)$ (and a dual homomorphism
    $c_0^t:Y\to A_0(k)$) defining a semiabelian variety $G_0$ (and a
    dual semiabelian variety $G^t$).
  \item a trivialization of the biextension $\tau_0:1_{X\times
      X}=1_{Y\times X}\to (c^t\times c)^*\cP_{A_0}^{-1}$, where
    $\cP_{A_0}$ is the Poincare bundle.
  \end{enumerate}

  When the abelian part $A_0$ is trivial, $\tau_0$ becomes simply a
  bilinear symmetric function $b_0:X\times X\to k^*$.
\end{saynum}

\begin{saynum}
  We now would like to explain how to associate this data to a curve
  $C$.  Part of this description can already be found in
  \cite[\S18]{Namikawa_NewComp2} and
  \cite[9.D]{Namikawa_ToroidalCompSiegel} where it is attributed to
  Mumford.

  \begin{enumerate}
  \item $X=H_1(\Gamma(C),\bZ)$. 
  \item By choosing arbitrarily an orientation on $\Gamma$, we get a
    natural embedding of $X$ in a free abelian group
    $C_1(\Gamma(C),\bZ)=\oplus\bZ e_j$, each $e_j$ corresponds to an
    edge of $\Gamma$. The form $B$ is the restriction to $H_1$ of the
    standard Euclidean form on $C_1$.
  \item an abelian variety $A_0$ is $\Pic^0(\wC)$. Instead of line
    bundle on $A_0$ we consider $B_0=Pic^{g-1}(\wC)$ and the natural
    line bundle $M_0$ on it defines by the theta divisor. A choice of
    an isomorphism $A_0\to B_0$ doesn't matter.
  \item every element of $H_1(\Gamma)$ defines a cycle of multidegree
    $(0,\dots,0)$ on $\wC$, i.e. an element of $A_0^t$. This gives the
    homomorphism $c$.
  \item Finally, the map $\tau$ is the most interesting part. The
    quick answer is that $\tau$ is given by a ``generalized
    crossratio''. 
  \end{enumerate}
\end{saynum}

\begin{saynum}
  Let $f,g$ be two meromorphic functions on a smooth projective curve
  $X$ with disjoint divisors. Then, defining
  \begin{displaymath}
    (f,g)=f(div\, g)= \prod_{x\in C} f(x)^{v_x(g)},
  \end{displaymath}
  one has $(f,g)=(g,f)$ according to A.Weil. For $f=(z-a)/(z-b)$,
  $g=(z-c)/(z-d)$ this is nothing but the usual crossratio.
  
  In \cite[XVII]{SGA4} Deligne showed that to arbitrary two invertible
  sheaves $L,M$ on $X$ and their meromorphic sections $f,g$ with
  disjoint divisors one can associate an element $(f,g)$ of a certain
  one-dimensional vector space $(L,M)$. These one-dimensional vector
  spaces are bilinear and symmetric in $L,M$ (in the case of degree 0
  they form a symmetric biextension of $\Pic^0\times\Pic^0$) and
  $(f,g)=(g,f)$ if $\deg M\cdot\deg L$ is even and $=-(g,f)$
  otherwise. We will call this pairing Deligne symbol. A very nice
  summary of its properties can be found in \cite{BeilinsonManin86}.
\end{saynum}

\begin{saynum}
  Now for every two distinct elements $e_k,e_l$ of the standard basis
  in $C_1(\Gamma)$ we have two divisors of the total degree 0 on
  $\wC=\cup \wC_i$. This defines a one-dimensional vector space
  $V_{k,l}$ and an element $(e_k,e_l)$ in it. If $k=l$, we still have
  a vector space $V_{k,k}$ but $(e_k,e_k)$ is undefined. {\em We
    define it arbitrarily.\/}
  
  In particular, restricting this to $X\times X\subset C_1\times C_1$,
  we obtain a pairing on $X\times X$ with the values in a certain
  collection of one-dimensional vector spaces. Because every element
  of $H_1$ has degree 0 on each irreducible component $\wC_i$, this
  pairing is symmetric.  It can be checked that these one-dimensional
  vector spaces are the fibers of $(c^t\times c)^*\cP_{A_0}^{-1}$,
  where $\cP_{A_0}$ is the Poincare bundle (= Weil biextension) on
  $A_0\times A_0$. This defines the trivialization $\tau_0$.

  In the case $A_0=0$, i.e. when all $C_i$ are rational, $\tau_0=b_0$
  is a product of the usual cross ratios.
\end{saynum}

\begin{saynum}\label{saynum:independence_of_the_choice}
  Our definition seemingly depends on a choice of $(e_k,e_k)$.
  However, by \cite{AlexeevNakamura96} an SQAV depends only on the
  residue class of $\tau_0$ modulo the following equivalence relation.
  $\tau_0$ can be replaced by
  \begin{displaymath}
    \tau'_0(x,y)=\tau_0(x,y)\cdot
    c^{B_1(x,y)}
  \end{displaymath}
  for any $c\in k$ and any symmetric positive bilinear form $B_1$
  defining the same Delaunay decomposition as $B$ (for the definitions
  of the Delaunay decompositions, see \cite{AlexeevNakamura96} or
  \cite{OdaSeshadri79}).
  
  The independence of the choice of $(e_k,e_k)$ now follows because on
  $C_1(\bR)$ the standard Euclidean form and the form $\sum
  \lambda_iz_i^2$ for any $\lambda_i>0$ define the same Delaunay
  decomposition. The Delaunay cells are the standard cubes and their
  faces.  Because, as one can easily show (\cite[3.2]{OdaSeshadri79}
  or \cite[\S18]{Namikawa_ToroidalCompSiegel}) $C_1(\Gamma(C),\bZ)\cap
  H_1(\Gamma(C),\bR)=H_1(\Gamma(C),\bZ)$, the Delaunay decomposition
  of $\XR=H_1(\bR)$ is the intersection of this standard Delaunay
  decomposition with $H_1(\bR)$.
  
  Therefore, every cell has the following simple description. For each
  $1\le i\le\dim C_1$ we choose an integer $n_i$ and two
  numbers $a_i,b_i$ with either $a_i=b_i=n_i$ or $a_i=n_i$ and
  $b_i=n_i+1$. Then we obtain the cell $\sigma$ in $C_1(\bR)$ defined
  by the inequalities
  \begin{displaymath}
    a_i \le z_i \le b_i 
  \end{displaymath}
  and the cell in $\sigma\cap H_1(\bR)$ in $H_1(\bR)$ (it may be
  empty).
\end{saynum}

\setcounter{example}{3}
\begin{example}
  In this case $H_1=C_1$ and we have the standard Euclidean space
  $\bR^g\supset\bZ^g$. The Delaunay cells are standard cubes and their
  faces. Modulo the translation by $\bZ^g$ there are exactly $\binom
  nk$ such cells of codimension $k$. These numbers are the same as in
  section \ref{sec:first description}.
  
  Further assume that the curve $C$ is rational for simplicity. Then
  the symmetric bilinear form $b_0$ is defined by $n(n-1)/2$
  crossratios $(e_i,e_j)$, $i< j$.
\end{example}

\setcounter{example}{4}
\begin{example}
  In this case $H_1(\Gamma(C),\bZ)\subset C_1(\Gamma(C),\bZ)$ is the
  hyperplane $\{x_1+x_2+x_3=0\}$. The Delaunay decomposition is the
  decomposition of $\bR^2$ into unilateral triangles. Modulo the
  translations there are two cells of dimension 2, 3 cells of
  dimension 1 and 1 cell of dimension 0. These numbers are the same as
  in section \ref{sec:first description}.
  
  This SQAV does not depend on the form $\tau_0$ as all the $3=\dim
  S^2(H_1)$ choices are killed by the $3=\dim C_1$ choices for
  $(e_k,e_k)$.
\end{example}

\setcounter{example}{5}
\begin{example}
  In this case $H_1(\Gamma(C),\bZ)\subset C_1(\Gamma(C),\bZ)$ is the
  hyperplane $\{x_1+...+x_n=0\}$. The lattice is the standard lattice
  $A_n$. It can be checked that the number of $k$-dimensional cells is
  given by the same formula as in the section \ref{sec:first
    description}.
\end{example}

\begin{saynum}
  \cite{AlexeevNakamura96} gives a stratification of an SQAV into
  locally closed subschemes which are homogeneous spaces over a
  semiabelian variety. A stratum of dimension $n$ corresponds to a
  Delaunay cell of dimension $n-\dim A_0$.
  
  On the other hand, in section \ref{sec:first description} we have
  given a similar description for $\Jac_{g-1}$ and the semiabelian
  variety is $\Pic^0 C$. In all the above examples the numbers of
  strata of each dimension in both descriptions are the same.  We now
  would like to relate the two descriptions explicitly.
\end{saynum}

\begin{saynum}
  Consider an arbitrary orientation of the generating subgraph
  $\Gamma'\subset\Gamma$. By \ref{defnprop} it corresponds to a
  semistable multidegree $\ud'$ of the graph $\Gamma(C)$. Now,
  depending on whether the edge $e_i$ is oriented the ``right'' way
  (the same that we used defining $H_1$), the ``wrong'' way, or not
  present at all, choose $a_i=0,b_i=1$, $a_i=-1,b_i=0$ or $a_i=b_i=0$.
  This gives a Delaunay cell $\sigma$ of $C_1(\bR)$ and the Delaunay
  cell $\sigma\cap H_1(\bR)$ of $H_1(\bR)$ as in
  \ref{saynum:independence_of_the_choice}.  Moreover, the orientation
  is stable iff
  \begin{displaymath}
    \dim\sigma=\dim\sigma\cap H_1(\bR)
  \end{displaymath}
\end{saynum}

\begin{proof}
  The above is Oda and Seshadri's description of the stratification of
  $\Jac_{\phi}$ for the case $\phi=\partial e(J)/2$, in which case the
  Namikawa-Delaunay decomposition of \cite{OdaSeshadri79} coincides
  with the Delaunay decomposition we have described above.  Therefore,
  everything follows from \cite{OdaSeshadri79} and the following lemma
\end{proof}

\begin{lem}
  $\Jac_{g-1}$ corresponds to the choice $\phi=\partial e(J)/2$ in the
  notations of \cite{OdaSeshadri79}. 
\end{lem}
\begin{proof}
  Follows directly from \ref{thm:Jacs_are_the_same}.
\end{proof}


\ifx\undefined\bysame
\newcommand{\bysame}{\leavevmode\hbox to3em{\hrulefill}\,}
\fi

\end{document}